\documentclass{amsart}
\usepackage{amsrefs, amsmath,amsthm,amssymb,enumerate,mathrsfs}
\usepackage{xcolor}
\usepackage[colorlinks=true,  linkcolor=red, citecolor=blue]{hyperref}
\newtheorem{theorem}{Theorem}[section]
\newtheorem{lemma}[theorem]{Lemma}

\newtheorem{remark}[theorem]{Remark}

\numberwithin{equation}{section}

\begin{document}
\title[]{Geometric inequalities for convex spacelike hypersurface in de Sitter space}

\author{Yandi Dong}
\address{School of Mathematics and Statistics,
Shandong University of Technology, Zibo 255000, China}
\email{ydong25@163.com}

\author{Kuicheng Ma}
\address{School of Mathematics and Statistics,
Shandong University of Technology, Zibo 255000, China}
\email{kchm@sdut.edu.cn}

\subjclass[2000]{53C44, 35K55}



\keywords{Locally constrained flows; Quermassintegrals; Alexandrov-Fenchel inequalities; Weighted curvature integrals; Maximum principle}

\begin{abstract}
In this paper, the long-time existence and convergence results are derived for locally constrained flows with initial value some compact spacelike hypersurface that is suitably pinched in the de Sitter space. As applications, geometric inequalities related to the quermassintegrals as well as the weighted curvature integrals are established.
\end{abstract}

\maketitle


\par
\section{Introduction}
Let $n$ be an integer larger than or equal to $2$, let $(a,b)\subset\mathbb{R}^+$ be some open interval and $\lambda: (a,b)\rightarrow \mathbb{R}^+$ be some smooth function such that $\lambda'>0$. Let $(\mathbb{S}^n, \sigma)$ be the standard sphere of dimension $n$ and $\epsilon$ be either $1$ or $-1$. Let $N$ be the product manifold
\begin{align}
(a,b)\times\mathbb{S}^n\nonumber
\end{align}
together with its structure
\begin{align}
\bar g= \epsilon \, dr^2+\lambda^2(r)\sigma,
\end{align}
which is usually referred to as a semi-Riemannian warped product. It is well known that the Euclidean space $\mathbb{R}^{n+1}$, the spherical space $\mathbb{S}^{n+1}$, the hyperbolic space $\mathbb{H}^{n+1}$ and the de Sitter space $\mathbb{S}^n_1$ can be recovered in this manner. Let $\bar\nabla$ be the Levi-Civita connection uniquely determined by $\bar g$. It is direct to check that the global vector field $V=\lambda\partial_r$ is conformal Killing. More precisely, there holds
\begin{align}\label{ck}
\bar\nabla V=\lambda'\mathcal{I},
\end{align}
where $\mathcal{I}$ is the identical transformation. From now on, $N$ is assumed further to be of constant sectional curvature $\bar K$.

Let $\Sigma$ be a compact spacelike hypersurface smoothly embedded in $N$ and $\nu$ be the particular unit normal along $\Sigma$. The generalized support function is defined as
\begin{align}
u=\epsilon ~\bar g (\lambda\partial_r,\nu),
\end{align}
which is positive due to the particular choice of the unit normal. Eigenvalues of the Weingarten transformation $\mathcal{W}=\bar\nabla\nu$, usually referred to as principal curvatures, are important geometric invariants.  For $k\in \{1,\cdots,n\}$, let $E_k$ be the $k$-th normalized elementary symmetric polynomial in the principal curvatures. For convenience, set $E_0=1$ and $E_{n+1}=0$.

A compact spacelike hypersurface in $N$ is said to be strictly convex if its principal curvatures are everywhere positive. More generally, a spacelike hypersurface in $N$ is $k$-convex if the $n$-tuple comprised of its principal curvatures always belongs to the cone
\begin{align}
\Gamma_k^+=\{(\kappa_1,\cdots,\kappa_n)~|~ E_m(\kappa_1,\cdots,\kappa_n) >0, ~m=1,\cdots,k\}\nonumber
\end{align}
which contains the positive cone
\begin{align}
\Gamma^+=\{(\kappa_1,\cdots,\kappa_n)~|~ \kappa_1>0,\cdots,\kappa_n>0\}.\nonumber
\end{align}
It is well known that $n$-convex is strictly convex, $1$-convex is usually referred to as mean convex. In particular, for hypersurfaces in hyperbolic space $\mathbb{H}^{n+1}$, both the horo-convexity and the static convexity are occasionally assumed.

Quermassintegrals of the closed spacelike hypersurface $\Sigma$ in $N$ are defined as
\begin{align}
\mathcal{A}_{-1}(\Sigma)&=(n+1)\int_{\hat\Sigma}~~ \text{dvol},~\;~\;~ \mathcal{A}_0(\Sigma)=\int_\Sigma ~~d\mu,\nonumber\\
\mathcal{A}_k(\Sigma)&=\int_\Sigma E_k ~d\mu+\epsilon \bar K \frac{k}{n+2-k}\mathcal{A}_{k-2}(\Sigma),~~k=1,\cdots, n,\nonumber
\end{align}
where $\hat\Sigma\subset N$ is the region enclosed by $\Sigma$ and $\{0\}\times \mathbb{S}^n$. Let
\begin{align}
\mathcal{B}_{-1}(\Sigma)&=(n+1)\int_{\hat\Sigma} \lambda' \text{dvol}+\lambda^{n+1}(0)\omega_n, \nonumber\\
\mathcal{B}_k(\Sigma)&=\int_\Sigma \lambda' E_k ~d\mu,~~\;~~ k=0,1, \cdots, n,\nonumber
\end{align}
which are collectively referred to as weighted curvature integrals. Here and in the sequel, $\omega_n$ stands for the area of the standard sphere $(\mathbb{S}^n, \sigma)$.

One current topic in differential geometry and geometric analysis is to establish optimal inequalities between quermassintegrals, also known as Alexandrov-Fenchel inequalities or isoperimetric type inequalities, using parabolic method. Consider the variation of some spacelike hypersurface in its normal direction, along which two quermassintegrals possess opposite monotonicity. Now it suffices to show such a variation will converge to some expected limit. Huisken \cite{GH} introduced a volume preserving mean curvature flow in $\mathbb{R}^{n+1}$ and his main result is an optimal inequality between quermassintegrals $\mathcal{A}_{-1}$ and $\mathcal{A}_0$ for strictly convex hypersurface there. See \cite{Mc1} and \cite{WX} for extensions. Instead of adding some globally constrained term, Guan, Li and Wang \cite{GLW} designed the locally constrained mean curvature flow
\begin{align}\label{flow1}
\frac{\partial}{\partial t} X=(\lambda'-uE_1)\nu,
\end{align}
in Riemannian warped products, which turns out to be a more powerful tool in the classical isoperimetric problem in the sense that it only requires the evolving hypersurface to be star-shaped. In fact, under mild assumptions over the ambient geometry, flow \eqref{flow1}, starting from any star-shaped hypersurface, exists eternally and converges to some coordinate slice. Moreover, along flow \eqref{flow1}, $\mathcal{A}_{-1}$ is preserved while the area $\mathcal{A}_0$ is monotone decreasing. Subsequently, Brendle, Guan and Li \cite{BGL} investigated the locally constrained inverse curvature flow
\begin{align}\label{flow2}
\frac{\partial }{\partial t} X=\bigg(\frac{\lambda'E_k}{E_{k+1}}-u\bigg)\nu,
\end{align}
in $\mathbb{H}^{n+1}$, aiming to weaken the convexity assumption imposed by Wang and Xia in \cite{WX}. Unfortunately, the uniform gradient estimate along flow \eqref{flow2} is still open in general cases, see \cite{GL} for details. As an application of flow \eqref{flow2}, Hu, Li and Wei \cite{HLW} provided a new proof for the main result shown in \cite{WX}. In fact, they discovered that the horo-convexity is preserved along the considered flow with the help of the maximum principle for tensor developed by Andrews \cite{BA}. Applications of locally constrained flows in establishing Alexandrov-Fenchel inequalities for strictly convex hypersurface in $\mathbb{S}^{n+1}$ can be traced in \cite{CGLS,CS}. More recently, locally constrained flows were successfully extended to isoperimetric type problems concerning the closed spacelike hypersurface in generalized Robertson-Walker spacetimes and that in $\mathbb{S}^n_1$ respectively by Lambert and Scheuer \cite{BJ} and Scheuer \cite{JS}. More precisely, as an application of the locally constrained mean curvature flow
\begin{align}
\frac{\partial }{\partial t} X=(uE_1-\lambda')\nu,\nonumber
\end{align}
Lambert and Scheuer \cite{BJ} established for any closed spacelike hypersurface $\Sigma$ in $\mathbb{S}^n_1$ the isoperimetric inequality
\begin{align}\label{isoi}
\mathcal{A}_0(\Sigma)\leq \varphi_0\circ\varphi_{-1}^{-1}(\mathcal{A}_{-1}(\Sigma)).
\end{align}
In this paper, we continue to study the locally constrained curvature flow
\begin{align}\label{FL-3}
\frac{\partial}{\partial t}X=\bigg(\frac{uE_k}{E_{k-1}}-\lambda'\bigg)\nu,
\end{align}
with $k\in\{1,\cdots,n\}$ and initial value a closed spacelike strictly convex hypersurface in de Sitter space $\mathbb{S}^n_1$. It turns out that an important feature of flow \eqref{FL-3} is the preservation of the strict convexity. But unfortunately, the long-time existence of flow \eqref{FL-3} is not available for $k\in\{2,\cdots, n\}$, even under the strict convexity assumption. As an application, we show

\begin{theorem}
Let $\Sigma_0$ be a closed spacelike hypersurface in de Sitter space $\mathbb{S}^n_1$, then
\begin{align}\label{DM}
\mathcal{A}_{-1}(\Sigma_0)\leq \varphi_{-1}\circ\psi_{-1}^{-1}(\mathcal{B}_{-1}(\Sigma_0)),
\end{align}
where $\varphi_{-1}$ and $\psi_{-1}$ are monotone increasing functions such that equality in \eqref{DM} holds if and only if $\Sigma_0$ is some coordinate slice. Assume further that $\Sigma_0$ is strictly convex, then there exists for each $l\in\{0,\cdots, n-2\}$ some monotone increasing function $\varphi_l$ such that
\begin{align}\label{DM-1}
\mathcal{A}_l(\Sigma_0)\leq \varphi_l\circ\varphi_{-1}^{-1}(\mathcal{A}_{-1}(\Sigma_0)),
\end{align}
and the equality holds if and only if $\Sigma_0$ is some coordinate slice.
\end{theorem}
\begin{remark}
Based on the dual relationship between closed spacelike strictly convex hypersurfaces in $\mathbb{S}^n_1$ and those in $\mathbb{H}^{n+1}$, together with their Blaschke-Santal\'o type inequalities, Hu and Li \cite{HL2} derived the Alexandrov-Fenchel inequalities
\begin{align}\label{la0}
\mathcal{A}_l\leq \varphi_l\circ\varphi_0^{-1}(\mathcal{A}_0), ~\;~l\in\{1,\cdots, n-1\}
\end{align}
for closed strictly convex spacelike hypersurfaces in $\mathbb{S}^n_1$. Comparing \eqref{isoi}, it follows that \eqref{la0} is better than \eqref{DM-1}. For $l=n-1$, \eqref{la0} was also addressed in \cite{M23}. On the other hand, for $l\in\{1,2\}$, \eqref{la0} holds actually for $l$-convex hypersurface, see \cite{M24-} and \cite{JS}.
\end{remark}
It turns out that locally constrained flows are as well powerful in establishing optimal inequalities between weighted curvature integrals. In \cite{SX}, Scheuer and Xia obtained the optimal inequality
\begin{align}
\int_\Sigma \lambda' E_1 d\mu-(n+1)\int_{\hat \Sigma} \lambda' \text{dvol} \geq \omega_n^{\frac{2}{n+1}}\bigg((n+1)\int_{\hat\Sigma} \lambda' \text{dvol}\bigg)^{\frac{n-1}{n+1}}\nonumber
\end{align}
for closed mean convex hypersurface in $\mathbb{H}^{n+1}$ by investigating the flow
\begin{align}\label{SC}
\frac{\partial}{\partial t} X=\bigg(\frac{1}{E_1}-\frac{u}{\lambda'}\bigg)\nu.
\end{align}
Recently, Hu and Li \cite{HL} established the optimal inequalities between $\mathcal{B}_{-1}(\Sigma)$ and $\mathcal{B}_k(\Sigma)$ for static convex hypersurface $\Sigma$ in $\mathbb{H}^{n+1}$, where $k\in\{0,\cdots,n\}$. In fact, they studied the locally constrained flow
\begin{align}
\frac{\partial }{\partial t} X=\bigg(1-\frac{uE_1}{\lambda'}\bigg)\nu\nonumber
\end{align}
and showed the preservation of the static convexity applying Andrews' maximum principle for tensor. In fact, Hu and Li addressed further the preservation of static convexity along variants of flow \eqref{SC}, replacing $E_1$ with $E_k/E_{k-1}$, and derived optimal inequalities between quermassintegrals and weighted curvature integrals. See \cite{HLW} and \cite{WZh} for similar results under horo-convexity assumption. Noticing that for hypersurfaces in $\mathbb{H}^{n+1}$, static convexity is stronger than mean convexity but weaker than horo-convexity.

Inspired by the idea in \cite{HL}, we shall assume the pinching condition
\begin{align}\label{PC}
\mathcal{W}\leq \Theta \mathcal{I}
\end{align}
for spacelike hypersurface in $\mathbb{S}^n_1$ and study in this paper the locally constrained inverse curvature flow
\begin{align}\label{fl1}
\frac{\partial}{\partial t} X=\bigg(1-\frac{E_{k-1}}{\Theta E_k}\bigg)\nu=\mathcal{F}\nu,
\end{align}
with initial value any closed spacelike $k$-convex hypersurface satisfying \eqref{PC}, and
\begin{align}
\Theta=\frac{u}{\lambda'}.\nonumber
\end{align}
Indeed, the pinching condition \eqref{PC} was introduced by the author in \cite{M24} to solve the weighted isoperimetric problem in $\mathbb{S}^n_1$.

\begin{theorem}
Let $X_0$ be a closed $k$-convex spacelike embedding into $\mathbb{S}^n_1$ such that on which \eqref{PC} holds, then with initial value $X_0$, flow \eqref{fl1} exists for all positive time and converges to some coordinate slice.
\end{theorem}

As an application, we obtain a family of optimal inequalities concerning the weighted curvature integrals.
\begin{theorem}
Let $\Sigma$ be a closed spacelike strictly convex hypersurface in $\mathbb{S}^n_1$ such that on which \eqref{PC} holds, then
\begin{align}\label{DM-2}
\mathcal{B}_n(\Sigma)\leq \psi_{n}\circ\psi_l^{-1}\big(\mathcal{B}_{l}(\Sigma)\big)
\end{align}
for each $l\in\{1,\cdots, n-1\}$ and
\begin{align}
\mathcal{A}_{n-1}(\Sigma) \leq \varphi_{n-1}\circ \psi^{-1}_{n}(\mathcal{B}_{n}(\Sigma)),\nonumber
\end{align}
where $\varphi_{n-1}$ and $\psi_1,\cdots,\psi_{n}$ are monotone functions such that the equality holds if and only if $\Sigma$ is some coordinate slice.
\end{theorem}
\begin{remark}
In fact, the Alexandrov-Fenchel type inequalities \eqref{DM-2} follow from those in \cite{HL} established for static convex hypersurface in $\mathbb{H}^{n+1}$, according to the dual relationship. In other words, we provide a new proof for the Alexandrov-Fenchel type inequalities in \cite{HL}.
\end{remark}

The rest of this paper is organized as follows. In Section \ref{S2}, we recall structural equations for spacelike hypersurface in the semi-Riemannian warped product and derive the uniform $C^0$ estimates; In Section \ref{S3}, we derive evolution equations along flow \eqref{fl1}  for important geometric quantities and obtain the gradient estimate by applying the classical maximum principle; In Section \ref{S4}, we show the preservation of $k$-convexity as well as the pinching condition \eqref{PC} along flow \eqref{fl1}; In Section \ref{S6}, we obtain the convergence of flow \eqref{fl1} and show the monotonicity of the weighted curvature integrals by assuming further that the spacelike hypersurface is strictly convex.

\section{Preliminaries}\label{S2}
\subsection{Semi-Riemannian geometry} For the semi-Riemannian manifold $N$, let
\begin{align}
\bar R(Y,Z)W =\bar\nabla_Y\bar\nabla_Z W-\bar\nabla_Z\bar\nabla_Y W-\bar\nabla_{[Y,Z]}W\nonumber
\end{align}
be its Riemann curvature tensor of type $(0,3)$ and then
\begin{align}
\bar R(U,W,Y,Z)=\bar g(U, R(Y,Z)W)\nonumber
\end{align}
Let $\Sigma$ be a spacelike hypersurface smoothly embedded in $N$, with the induced metric $g$. The vector-valued second fundamental form $B$ can be defined via the following Gauss formula
\begin{align}
\bar\nabla_Y Z = \nabla_YZ+B(Y,Z),\nonumber
\end{align}
where $\nabla$ is the Levi-Civita connection uniquely determined by $g$. Here and in the sequel, we will not distinguish a vector field from its push-forward. Let
\begin{align}
B(Y,Z) =-\epsilon h(Y,Z)\nu,\nonumber
\end{align}
then
\begin{align}
h(Y,Z) =g(\mathcal{W}(Y),Z) =g(Y,\mathcal{W}(Z)).\nonumber
\end{align}
As a consequence,
\begin{align}
\bar R(Y, Z)U&=R(Y, Z)U-\epsilon h(Z, U)\mathcal{W}(Y)+ \epsilon h(Y,U)\mathcal{W}(Z)\nonumber\\
&-\epsilon(\nabla_Yh)(Z,U)\nu +\epsilon ( \nabla_Zh)(Y,U)\nu,\nonumber
\end{align}
from which it yields the Gauss equation
\begin{align}
R(W,U,Y,Z)=R(W,U,Y,Z)-\epsilon h(Z,U)h(Y,W)+\epsilon h(Y,U)h(Z,W)
\end{align}
and the Codazzi equation
\begin{align}
(\nabla_Y h)(Z,U)-(\nabla_Z h)(Y,U)= -\bar R(\nu, U, Y, Z).
\end{align}
For symmetric $(0,2)$ tensor $T$, the Ricci identity addresses the rule for exchanging the order in taking second-order covariant derivative. More precisely, there holds
\begin{align}
(\nabla^2 T)(Y,Z;U,W)-(\nabla^2 T)(Y,Z;W,U)=T(R(U,W)Y, Z)+T(Y, R(U,W)Z).\nonumber
\end{align}
It is well known that the de Sitter space $\mathbb{S}^n_1$ can be characterized as the Lorentzian warped product
\begin{align}
(0,+\infty ) \times \mathbb{S}^n, ~\;~\;~\bar g= dr^2+ \lambda^2(r)\sigma,
\end{align}
where $\lambda(r)=\cosh(r)$. It turns out that $\mathbb{S}^n_1$ is of constant sectional curvature $1$, or equivalently
\begin{align}
\bar R(W,U,Y,Z)=g(Y,W)g(Z,U)-g(Z,W)g(Y,U).\nonumber
\end{align}
\subsection{Elementary symmetric functions and related properties} From
\begin{align}\label{deco}
E_{k-1}=\frac{n-k+1}{k}\frac{\partial E_k}{\partial \kappa_i}+\kappa_i\frac{\partial E_{k-1}}{\partial \kappa_i}~\;~\;~ (k=1,\cdots,n+1)
\end{align}
and the homogeneity, it yields
\begin{lemma}
For $(\kappa_1,\cdots, \kappa_n)\in\mathbb{R}^n$,
\begin{align}
\sum^n_{i=1}\frac{\partial F }{\partial \kappa_i}&=k-(k-1)\frac{E_k E_{k-2}}{E^2_{k-1}},\nonumber\\
\sum^n_{i=1}\kappa^2_i \frac{\partial F}{\partial \kappa_i}&=(n-k+1)F^2-(n-k)\frac{E_{k+1}}{E_{k-1}},\nonumber
\end{align}
here and in the sequel, $k\in\{1,\cdots, n\}$ and
\begin{align}
F=\frac{E_k}{E_{k-1}}.\nonumber
\end{align}
\end{lemma}
The well known Newton inequalities are formulated as
\begin{lemma}
For $(\kappa_1,\cdots,\kappa_n)\in\Gamma_k$, there holds the inequalities
\begin{align}
\frac{E_{k+1}}{E_k}\leq \frac{E_k}{E_{k-1}}\leq \cdots\leq E_1,\nonumber
\end{align}
and the equality holds if and only if $\kappa_1=\cdots =\kappa_n$.
\end{lemma}
As a consequence, for $k\in\{1,\cdots, n\}$ and $(\kappa_1,\cdots,\kappa_n)\in\Gamma_k$,
\begin{align}\label{ki}
\sum^n_{i=1}\frac{\partial F }{\partial\kappa_i}\geq 1, ~\;~\sum^n_{i=1}\kappa^2_i \frac{\partial F}{\partial \kappa_i}\geq F^2.
\end{align}
Sometimes, it is more convenient to consider $F$ as function in the Weingarten transformation $\mathcal{W}$ and let
\begin{align}
F^{pq,rs}=g^{ls}g^{jq}\frac{\partial^2F}{\partial h^j_p \partial h^l_r},\nonumber
\end{align}
then there holds the identity
\begin{align}
F^{pq,rs}h_{pq;i}h_{rs;i}=\sum^n_{p, r=1}\frac{\partial^2F}{\partial\kappa_p\partial\kappa_r}h_{pp;i}h_{rr;i}+2\sum_{p>r} \frac{\frac{\partial F}{\partial \kappa_p}-\frac{\partial F}{\partial \kappa_r}}{\kappa_p-\kappa_r} h^2_{pr;i},\nonumber
\end{align}
where $h_{pq;i}$ be local components of $\nabla h$.

At the end of this subsection, we intend to emphasize that $F$ is strictly monotone increasing and concave in $\Gamma_k$, which turns to be crucial in the behaviors of flow \eqref{FL-3} and flow \eqref{fl1}.
\subsection{$C^0$ estimate} Since every compact spacelike hypersurface embedded in $\mathbb{S}^n_1$ can be written as graph, we may assume the evolving hypersurface along flow \eqref{fl1} as
\begin{align}
X=\big\{(r, \xi)~~|~~ \xi\in\mathbb{S}^n\big\},\nonumber
\end{align}
where $r$ is a family of smooth functions defined on $\mathbb{S}^n$, and it is spacelike in case the symmetric matrix
\begin{align}\label{im}
g_{ij}=\lambda^2\sigma_{ij}-r_i r_j
\end{align}
is positive definite, or equivalently,
\begin{align}
1-\lambda^{-2}|Dr|^2>0,\nonumber
\end{align}
where $D$ is the Levi-Civita connections determined by  $\sigma$ and $|Dr|^2$ is the squared norm with respect to the spherical metric $\sigma$.

Let
\begin{align}\label{nv}
\nu=\frac{1}{\upsilon}\big(\partial_r+\lambda^{-2} Dr\big),
\end{align}
be the future-directed timelike unit normal along spacelike hypersurface $X$, where
\begin{align}
\upsilon=\sqrt{1-\lambda^{-2}|Dr|^2}.\nonumber
\end{align}
And the generalized support function is defined as
\begin{align}
u=-\bar g (\lambda\partial_r,\nu)=\frac{\lambda}{\upsilon}.\nonumber
\end{align}

The second fundamental form, with respect to the normal in \eqref{nv}, determines the symmetric matrix
\begin{align}\label{sff}
h_{ij}=\frac{1}{\upsilon}\big(r_{,ij}+\lambda\lambda'\sigma_{ij}-2\lambda^{-1}\lambda'r_ir_j\big),
\end{align}
where the covariant derivative $r_{,ij}$  is with respect to $\sigma$. Furthermore,
\begin{align}\label{Hor}
h_{ij}=\upsilon\big(r_{;ij}+\lambda^{-1}\lambda'g_{ij}+\lambda^{-1}\lambda'r_ir_j\big),
\end{align}
where the covariant derivative $r_{;ij}$ is with respect to the induced metric \eqref{im} and
\begin{align}
g^{ij}=\lambda^{-2}\bigg(\sigma^{ij}+\frac{r^i r^j}{\lambda^2\upsilon^2}\bigg).\nonumber
\end{align}

Applying De Turck's trick, flow \eqref{fl1} is equivalent to
\begin{align}
\frac{\partial }{\partial t} r=\upsilon \bigg(1-\frac{1}{\Theta F}\bigg),\nonumber
\end{align}
and at the point where $r$ attains its spatial maximum, $\upsilon=1$ by definition, principal curvatures are no larger than $\lambda^{-1}\lambda'$ due to \eqref{sff} and the classical maximum principle while the generalized support function is equal to $\lambda$ there. Hence, $\Theta F$ is no larger than $1$, where the monotonicity and homogeneity of function $F$ are involved. It follows immediately that the spatial maximum of function $r$ is decreasing along flow \eqref{fl1}. Similarly, we can conclude that the spatial minimum of function $r$ is increasing. Hence
\begin{align}
\min_{\mathbb{S}^n} r(\xi ,0) \leq r(\xi, t)\leq \max_{\mathbb{S}^n}  r(\xi, 0).\nonumber
\end{align}
Noticing that the $C^0$ a priori estimate for flow \eqref{FL-3} follows from the same argument.
\section{evolution equations and gradient estimate}\label{S3}
From \eqref{Hor}, it follows that
\begin{align}\label{How}
(\lambda')_{;pq}=uh_{pq}-\lambda'g_{pq}.
\end{align}
For $l\in\{1,\cdots,n\}$, contracting \eqref{How} with $E_l^{pq}$ and then taking integral over the closed spacelike hypersurface $\Sigma$ in $N$, it yields
\begin{align}\label{MI}
\int_{\Sigma} (uE_l-\lambda' E_{l-1})d\mu=0,
\end{align}
where \eqref{deco} and the homogeneity of the elementary symmetric function are used.

On the other hand, we have
\begin{align}\label{gos}
\nabla u=\mathcal{W}(\nabla\lambda')
\end{align}
and
\begin{align}\label{Hos}
u_{;pq}=-\lambda'h_{pq}+g(\nabla h_{pq}, \nabla\lambda')+uh_p^m h_{mq}.
\end{align}
From \eqref{gos}, it yields that
\begin{align}\label{got}
\nabla \Theta =(\mathcal{W} -\Theta\mathcal{I})(\nabla\ln \lambda').
\end{align}
Combining \eqref{How} and \eqref{Hos}, it turns out
\begin{align}
\nabla^2\Theta&=\Theta g-(1+\Theta^2)h+g(\nabla\ln\lambda', \nabla h)+\Theta h^2-\nabla\ln\lambda'\otimes \nabla \Theta,
\end{align}
where
\begin{align}
(\nabla\ln\lambda'\otimes \nabla \Theta)(Y,Z)=g(\nabla\ln\lambda',Y) g(\nabla \Theta, Z)+g(\nabla\ln\lambda', Z) g(\nabla\Theta,Y).\nonumber
\end{align}

Along flow \eqref{fl1}, we introduce the second order partial differential operator
\begin{align}
\mathcal{P}=\frac{\partial}{\partial t} -\frac{1}{\Theta F^2}F^{pq}\nabla_p\nabla_q -\frac{1}{\Theta^2F} g(\nabla \ln\lambda' , \nabla  ).\nonumber
\end{align}
From $\Theta>1$, together with the positivity and the monotonicity of $F$ in $\Gamma_k$, we can conclude that the second order operators $\mathcal{P}$ is parabolic.
\begin{lemma}\label{eow}
Along flow \eqref{fl1}, there holds the evolution equations
\begin{align}
\mathcal{P} \lambda'=u\bigg(1-\frac{2}{\Theta F}+\frac{1}{(\Theta F)^2} F^{pq} g_{pq}\bigg)-\frac{\lambda'}{\Theta^2F}||\nabla\ln\lambda'||^2,\nonumber
\end{align}
where the norm $||\cdot||$ is with respect to the induced metric on evolving hypersurface.
\end{lemma}
\begin{proof}
From \eqref{How} and
\begin{align}
\frac{\partial}{\partial t}\lambda'=-\bar g (\lambda\partial_r,\nu)\mathcal{F}=u\mathcal{F},\nonumber
\end{align}
the desired evolution equation follows.
\end{proof}

\begin{lemma}\label{eos}
Along flow \eqref{fl1}, there holds the evolution equation
\begin{align}
\mathcal{P} u=\lambda'\bigg(1-\frac{1}{F^2}F^{pq}h_p^m h_{mq}\bigg)-\frac{\lambda'}{\Theta F}||\nabla\ln\lambda'||^2.\nonumber
\end{align}
\end{lemma}
\begin{proof}
From \eqref{Hos}, together with either
\begin{align}
\frac{\partial }{\partial t} \nu=\nabla\mathcal{F}\nonumber
\end{align}
or further
\begin{align}
\frac{\partial }{\partial t} u=-\bar g(\lambda'\mathcal{F}\nu ,\nu)-\bar g(\lambda\partial_r, \nabla\mathcal{F})=\lambda'\mathcal{F}+g(\nabla\lambda',\nabla\mathcal{F}),\nonumber
\end{align}
the desired evolution equation follows.
\end{proof}
Applying the maximum principle to the evolution equation in Lemma \ref{eos}, we may conclude that the spatial maximum of the generalized support function is monotone decreasing. As a consequence, flow \eqref{fl1} preserves to be spacelike if it initially is.

Furthermore, combining Lemmata \ref{eow} and \ref{eos}, we obtain
\begin{lemma}\label{eot}
Along flow \eqref{fl1}, there holds the evolution equation
\begin{align}
\mathcal{P}\Theta=-\bigg(\Theta^2-\frac{2\Theta}{F}+\frac{1}{F^2}F^{pq}g_{pq}\bigg)+\bigg(1-\frac{1}{F^2}F^{pq}h_p^mh_{mq}\bigg)+\frac{2}{ F^2}F^{pq}(\ln \Theta)_{;p}(\ln \lambda')_{;q}.\nonumber
\end{align}
\end{lemma}

\begin{lemma}\label{eve1}
Let
\begin{align}
\mathcal{Q}=\frac{\partial }{\partial t}-uF^{pq}\nabla_p\nabla_q -g(F\nabla\lambda', \nabla),\nonumber
\end{align}
then there holds the evolution equation
\begin{align}
\mathcal{Q} u =-(\lambda')^2-||\nabla \lambda'||^2+2\lambda' uF-u^2 F^{pq}h_p^m h_{mq}\nonumber
\end{align}
along flow \eqref{FL-3}.
\end{lemma}
Similarly, the generalized support function is uniformly bounded along flow \eqref{FL-3}. For simplicity, the proof of Lemma \ref{eve1}, as well as those of Lemmata \ref{eve2} and \ref{eve3}, is omitted since they make no difference with those of Lemmata \ref{eow}, \ref{eos}, \ref{eosf} and \ref{eoc}, except replacing $\mathcal{F}$ therein with the speed of flow \eqref{FL-3}.
\section{preservation of the $k$-convexity and the pinching condition \eqref{PC} along flow \eqref{fl1}}\label{S4}
In this section, we show preservation of $k$-convexity along flow \eqref{fl1} first.

\begin{lemma}\label{eosf}
Along flow \eqref{fl1}, there holds the evolution equation
\begin{align}
\mathcal{P} h_{ij}&=\bigg(1+\frac{1}{\Theta F}\bigg)h_i^l h_{lj}-\bigg(\frac{1}{F}+\frac{1}{\Theta^2 F} +\frac{1}{\Theta F^2}F^{pq}h_p^m h_{mq}+\frac{1}{\Theta F^2}F^{pq}g_{pq}\bigg)h_{ij}\nonumber\\
&+\bigg(1+\frac{1}{\Theta F}\bigg)g_{ij}+\frac{1}{\Theta F^2}F^{pq,rs} h_{pq;i}h_{rs;j}-\frac{2}{\Theta F} (\ln F)_{;i}(\ln F)_{;j}\nonumber\\
&-\frac{1}{\Theta F}\big[(\ln \Theta)_{;i} (\ln F)_{;j}+ (\ln F)_{;i} (\ln \Theta)_{;j}\big]-\frac{2}{\Theta F}(\ln\Theta)_{;i}(\ln\Theta)_{;j}\nonumber\\
&-\frac{1}{\Theta F} \big[(\ln \lambda')_{;i}(\ln \Theta)_{;j}+(\ln \Theta)_{;i}(\ln \lambda')_{;j}\big].\nonumber
\end{align}
\end{lemma}

\begin{proof}
Taking derivative of the second fundamental form with respect to time $t$, we have
\begin{align}
\frac{\partial}{\partial t} h_{ij}=\mathcal{F}g_{ij}+\mathcal{F}_{;ij}+\mathcal{F}h_i^l h_{lj}.\nonumber
\end{align}
Using \eqref{How}, \eqref{Hos} and the Gauss equations, the Codazzi equations as well as the Ricci identities, the second order covariant derivative $\mathcal{F}_{;ij}$ can be expanded in detail. More precisely,
\begin{align}
\mathcal{F}_{;ij}&=\frac{1}{\Theta F^2} F^{pq} h_{ij;pq}+\frac{1}{\Theta^2 F} g (\nabla\ln \lambda',\nabla h_{ij}) + \frac{2}{\Theta F}h_i^l h_{lj}+\frac{2}{\Theta F}g_{ij}\nonumber\\
&-\frac{1+\Theta^2}{\Theta^2 F}h_{ij}-\frac{1}{\Theta F^2} \big(F^{pq}h_p^m h_{mq}+F^{pq}g_{pq}\big)h_{ij}+\frac{1}{\Theta F^2}F^{pq,rs} h_{pq;i}h_{rs;j}\nonumber\\
&-\frac{1}{\Theta F^2}\big[F_{;i}(\ln \Theta)_{;j}+(\ln \Theta)_{;i}F_{;j}\big]-\frac{2}{\Theta F^3} F_{;i}F_{;j}-\frac{2}{\Theta F} (\ln \Theta)_{;i} (\ln \Theta)_{;j}\nonumber\\
&-\frac{1}{\Theta F} \big[(\ln \Theta)_{;i}(\ln \lambda')_{;j}+(\ln \lambda')_{;i}(\ln \Theta)_{;j}\big].\nonumber
\end{align}
\end{proof}

Notice that the covariant derivatives of the Weingarten transformation is equal to that of the second fundamental form while
\begin{align}
\frac{\partial}{\partial t} h^i_j = -2\mathcal{F}h^i_lh^l_j +g^{il}\frac{\partial}{\partial t }h_{lj},\nonumber
\end{align}
from which it follows
\begin{lemma}\label{eoc}
Along flow \eqref{fl1}, there holds the evolution equation
\begin{align}
\mathcal{P} F&=(F^{pq}g_{pq}-1)-\Theta^{-2}+\bigg(\frac{2}{\Theta F}-1\bigg)F^{pq}h_p^m h_{mq}-\frac{2}{\Theta F}F^{pq}(\ln \Theta)_{;p} (\ln \lambda')_{;q}\nonumber\\
&-\frac{2}{\Theta F}F^{pq}(\ln \Theta)_{;p} (\ln \Theta)_{;q}-\frac{2}{\Theta F}F^{pq}(\ln \Theta)_{;p} (\ln F)_{;q}-\frac{2}{\Theta F}F^{pq}(\ln F)_{;p}(\ln F)_{;q}.\nonumber
\end{align}
\end{lemma}
\begin{proof}
Consider $F$ as function in the Weingarten transformation, then from the chain rule it yields that
\begin{align}
\frac{\partial }{\partial t} F = -2\mathcal{F}F^{pq}h_p^mh_{mq} +F^{ij}\frac{\partial}{\partial  t} h_{ij}.\nonumber
\end{align}
Substituting the time-derivative of the second fundamental form into the equation above, the desired evolution equation follows.
\end{proof}

Now we turn to the preservation of $k$-convexity. To this end, we consider the auxiliary function
\begin{align}
\omega=-\ln F-\ln \Theta-\ln (\lambda'-\delta),\nonumber
\end{align}
where
\begin{align}
\delta=\frac{1}{2}\min \lambda'.\nonumber
\end{align}
Without loss of generality, we may assume that the maximum of the auxiliary function is unbounded, or equivalently the function $F$ may approach to zero, since otherwise the proof is done. From Lemmata \ref{eot} and \ref{eoc}, it yields directly the evolution equation of this auxiliary function. Noticing that at the point where the auxiliary function attains its maximum,
\begin{align}
\nabla\ln F+\nabla\ln\Theta =-\nabla\ln (\lambda'-\delta),\nonumber
\end{align}
and as a consequence, there holds the inequality
\begin{align}
0&\leq -\frac{\delta}{\lambda'-\delta}\Theta-\Theta^{-1}-\frac{1}{F}(F^{pq}g_{pq}-1) -\frac{1}{F}\bigg(\frac{1}{\Theta F}-1\bigg) F^{pq}h^m_ph_{mq}\nonumber\\
&-\frac{\delta}{\lambda'-\delta }\frac{1}{\Theta F^2}F^{pq}g_{pq}+\frac{1}{\Theta^2F}+\frac{\lambda'}{\lambda'-\delta} \frac{1}{\Theta^2F}||\nabla \ln\lambda'||^2+\frac{2\delta}{\lambda'-\delta}\frac{1}{F}\nonumber
\end{align}
at the considered point and a contradiction occurs immediately.

For convenience, the maximum principle for tensors developed by Andrews in \cite{BA} is formulated here as
\begin{theorem}\label{MPT}
Let $S_{ij}$ be a smooth time-varying symmetric tensor field on a closed manifold $M$ satisfying
\begin{align}
\frac{\partial}{\partial t} S_{ij}=a^{kl}\nabla_k\nabla_l S_{ij} + b^k\nabla_k S_{ij}+ N_{ij},\nonumber
\end{align}
where $a^{kl}$ and $b^k$ are smooth, $\nabla$ is a smooth symmetric connection, and $a^{kl}$ is positive definite everywhere. Suppose
\begin{align}\label{NPC}
N_{ij}\eta^i\eta^j+\sup_{\Lambda} 2a^{kl}\big(2\Lambda_k^p\nabla_l S_{ip}\eta^i-\Lambda_k^p\Lambda_l^qS_{pq}\big)\geq 0,
\end{align}
whenever $S_{ij}\geq 0$ and $S_{ij}\eta^j=0$ and $\Lambda$ is an $n\times n$ matrix. If $S_{ij}$ is positive definite everywhere on $M$ at $t=0$, then it is positive on $M\times[0,T]$.
\end{theorem}

Let $S_{ij}=\Theta g_{ij}- h_{ij}$. From lemmata \ref{eot} and \ref{eosf}, it yields that
\begin{align}
\mathcal{P} S_{ij}&=-\bigg(1+\frac{1}{\Theta F}\bigg)S_i^l S_{lj}+\bigg(\frac{3}{F}-\frac{1}{\Theta^2 F}-\frac{F^{pq}h_p^m h_{mq}}{\Theta F^2}-\frac{F^{pq}g_{pq}}{\Theta F^2}\bigg) S_{ij}\nonumber\\
&+\frac{1}{\Theta F} \big[(\ln F+\ln u)_{;i}(\ln \Theta)_{;j}+(\ln F+\ln u)_{;j}(\ln\Theta)_{;i}\big]\nonumber\\
&+\frac{2}{\Theta F}(\ln F)_{;i}(\ln F)_{;j}+\frac{2}{F^2}F^{pq}(\ln \Theta)_{;p}(\ln \lambda')_{;q} g_{ij}-\frac{1}{\Theta F^2}F^{pq,rs} h_{pq;i}h_{rs;j},\nonumber
\end{align}
since
\begin{align}
\frac{\partial}{\partial t} g_{ij}=2\mathcal{F} h_{ij}\nonumber
\end{align}
and
\begin{align}
h_i^l h_{lj}=S_i^l S_{lj}-2\Theta S_{ij}+\Theta^2 g_{ij}.\nonumber
\end{align}

At the point where $S_{ij}$ develops a null vector $\eta$, i.e., $S_{ij}\eta^j=0$, we may choose the normal coordinates such that $g_{ij}=\delta_{ij}$. By continuity, we may assume the principal curvatures are mutually distinct and in decreasing order, that is $\kappa_1>\kappa_2>\cdots>\kappa_n$. The null vector condition $S_{ij}\eta^j=0$ now implies that $\eta=e_1$ and $S_{11}=\Theta-\kappa_1=0$ at the considered point. The terms involving $S_{ij}$ and $(S^2)_{ij}$ satisfy the null vector condition and can be ignored. Moreover,
\begin{align}
\nabla_1 (\ln  \Theta) =\frac{\kappa_1-\Theta}{\Theta}\nabla_1(\ln \lambda')=0.\nonumber
\end{align}

It turns out that one obstruction in applying Hamilton's maximum principle for tensor (see \cite{Ha}) is the term
\begin{align}
\frac{2}{F^2}F^{pq}(\ln \Theta)_{;p}(\ln \lambda')_{;q}=\frac{2}{F^2}\sum_{p=2}^nF^{pp}(\ln \Theta)_{;p}(\ln \lambda')_{;p}.\nonumber
\end{align}
Since $S_{11}=0$ and $\nabla S_{11}=0$, there holds
\begin{align}
&\frac{2}{\Theta F^2} F^{ij}\big( 2 \Lambda_i^p \nabla_j S_{1p} -\Lambda_i^p\Lambda_j^q S_{pq}\big)=\frac{2}{\Theta F^2}\sum_{i=1}^n\sum_{p=2}^n F^{ii}\big[2\Lambda_i^p \nabla_i S_{1p} -(\Lambda_i^p)^2 S_{pp}\big]\nonumber\\
&=\frac{2}{\Theta F^2}\sum_{i=1}^n\sum_{p=2}^n F^{ii}\bigg[\frac{\big(\nabla_i S_{1p}\big)^2}{S_{pp}} -\bigg(\Lambda_i^p-\frac{\nabla_i S_{1p}}{S_{pp}}\bigg)^2 S_{pp}\bigg].\nonumber
\end{align}
As a consequence, the supremum in \eqref{NPC} is now equal to
\begin{align}
\frac{2}{\Theta F^2}\sum_{i=1}^n\sum_{p=2}^n F^{ii}\frac{\big(\nabla_i h_{1p}\big)^2}{S_{pp}}.\nonumber
\end{align}
On the other hand, it yields from the concavity of $F$ further that
\begin{align}
&-\frac{1}{\Theta F^2}F^{pq,rs} h_{pq;1}h_{rs;1}+\frac{2}{\Theta F^2}\sum_{i=1}^n\sum_{p=2}^n F^{ii}\frac{\big(\nabla_i h_{1p}\big)^2}{S_{pp}}\nonumber\\
&\geq -\frac{2}{\Theta F^2}\sum_{p=2}^n\frac{F^{pp}-F^{11}}{\kappa_p-\kappa_1}(\nabla_p h_{11})^2+\frac{2}{\Theta F^2}\sum_{p=2}^n F^{11}\frac{\big(\nabla_p h_{11}\big)^2}{S_{pp}}\nonumber\\
&=\frac{2}{\Theta F^2}\sum_{p=2}^n \frac{F^{pp}\Theta^2_{;p}}{\Theta-\kappa_p},\nonumber
\end{align}
where $S_{11}=0$ and $\nabla S_{11}=0$ are used again.
Finally, it suffices to show
\begin{align}
\frac{2}{F^2}\sum_{p=2}^nF^{pp}(\ln \Theta)_{;p}(\ln \lambda')_{;p}+\frac{2}{F^2}\sum_{p=2}^n\frac{F^{pp}}{\Theta-\kappa_p}(\ln\Theta)_{;p}\Theta_{;p}\geq 0,\nonumber
\end{align}
which is trivial since $\Theta_{;p}=(\kappa_p-\Theta)(\ln\lambda')_{;p}$.
\section{Preservation of the strict convexity along flow \eqref{FL-3}}
\begin{lemma}\label{eve2}
Along flow \eqref{FL-3}, there hold the evolution equations
\begin{align}
\mathcal{Q} \lambda'=\lambda' u(F^{pq} g_{pq}-1)-F||\nabla\lambda'||^2\nonumber
\end{align}
and
\begin{align}
\mathcal{Q}F=uF \big(F^{pq}g_{pq}-1\big)+\lambda' \big(F^{pq}h_p^m h_{mq}-F^2\big)+2F^{pq}u_{;p}F_{;q}.\nonumber
\end{align}
\end{lemma}
It follows from the maximum principle that $F$ is uniformly bounded from below along flow \eqref{FL-3}. To derive the uniform upper bound for $F$, we consider the evolution equation of the auxiliary function
\begin{align}
\omega=\ln F+\ln u-\ln\lambda'.\nonumber
\end{align}
In fact, along flow \eqref{FL-3}, there holds the inequality
\begin{align}
0\leq-\bigg(u-\frac{\lambda'}{F}\bigg)F^{pq}h_p^m h_{mq}+\frac{u^2-1}{\lambda'}F-u+\frac{1}{u}\nonumber
\end{align}
at the global maximum point of $\omega$, where the critical point condition
\begin{align}
0=\nabla\ln F+\nabla\ln u-\nabla\ln \lambda'\nonumber
\end{align}
is used. Without loss of generality, we may assume the curvature function $F$ is unbounded at the considered point. But now contradiction occurs due to \eqref{ki}.

\begin{lemma}\label{eve3}
Along flow \eqref{FL-3}, there holds the evolution equation
\begin{align}
\mathcal{Q} h^i_j&=(\lambda'+uF)h^i_l h^l_j-(u+\lambda' F+u F^{pq}g_{pq}+u F^{pq}h_p^m h_{mq}) h^i_j\nonumber\\
&+2uF\delta^i_j+u_{;} ^{\; i} F_{;j}+F_{;} ^{\; i} u_{;j}+uF^{pq,rs} h_{pq;}^{\; \; \; \; i}h_{rs; j},\nonumber
\end{align}
\end{lemma}

Let $(b^i_j)$ be the inverse of $(h^i_j)$ such that $b^i_l h^l_j=\delta^i_j$. As a consequence, there holds
\begin{align}
\mathcal{Q} b^l_m&=-(\lambda'+uF)\delta^l_m+(u+\lambda' F+u F^{pq}g_{pq}+u F^{pq}h_p^m h_{mq})b^l_m-2uFb^l_i b^i_m\nonumber\\
&-b^l_i u_{;}^{\; i} b^j_m F_{;j}-b^l_i F_{;}^{ \; i} b^j_m u_{;j}-u \big(F^{pq,rs}+2F^{qs} b^{pr}\big) h_{pq;}^{\; \; \; \; i} b^l_i h_{rs; j} b^j_m\nonumber
\end{align}
along flow \eqref{FL-3}. Suppose that the convexity is lost somewhere, then the largest eigenvalue of $(b^i_j)$ blows up. But the largest eigenvalue of $(b^i_j)$ is in generally not smooth, we consider the auxiliary function
\begin{align}
\phi=\sup \{b_{ij} \xi^i\xi^j~|~ g_{ij}\xi^i\xi^j=1\}\nonumber
\end{align}
instead. At the point where $\phi$ attains its maximum, we choose local coordinates such that
\begin{align}
g_{ij}=\delta_{ij},~\;~ b_{ij}=\kappa_i^{-1}\delta_{ij}, ~\;~ \kappa_1\leq \cdots\leq \kappa_n.\nonumber
\end{align}
Around the considered point, fix $\xi=(1,0,\cdots,0)\in\mathbb{R}^n$ and let
\begin{align}
\omega=\frac{b_{ij}\xi^i\xi^j}{g_{ij}\xi^i\xi^j}.\nonumber
\end{align}
It turns out that at the considered point the auxiliary function $\omega$ satisfies the same evolution equation as $b^1_1$ and there holds the inequality
\begin{align}
0&\leq -(\lambda'+uF)+(u+\lambda' F+u F^{pq}g_{pq}+u F^{pq}h_p^m h_{mq})b^1_1-2uF(b^1_1)^2+\frac{u^2-\lambda^2}{2u} F,\nonumber
\end{align}
where the inverse concavity of the curvature function and \eqref{gos} are used. In fact,
\begin{align}
\big(F^{pq,rs}+2F^{qs} b^{pr}\big) h_{pq;}^{\; \; \; \; i} b^l_i h_{rs; j} b^j_m \geq \frac{2}{F} b^l_i F_{;}^{ \; i} F_{;j} b^j_m.\nonumber
\end{align}
Based on the a priori estimates obtained in previous sections, contradiction occurs if $b^1_1$ is permitted to be unbounded.
\section{convergences and geometric inequalities} \label{S6}
Now the mean convexity together with the pinching condition \eqref{PC} implies the uniform curvature estimate along flow \eqref{fl1}, from which the long-time existence follows, according to the regularity theory in \cite{NK}. In the particular case where $k=1$, the long-time convergence and convergence of flow \eqref{FL-3} was claimed in \cite{BJ}. In the section, we assume $k=1$ in \eqref{FL-3}.

According to the a priori estimates, the weighted curvature integral $\mathcal{B}_{k-1}$ is bounded along flow \eqref{fl1}. On the other hand, along flow \eqref{fl1},
\begin{align}
&\frac{d}{dt}\mathcal{B}_{k-1}(\Sigma_t)=k\int_{\Sigma_t} u E_{k-1}\mathcal{F} d\mu+(n-k+1)\int_{\Sigma_t}\lambda' E_{k}\mathcal{F}~d\mu\nonumber\\
&=k\int_{\Sigma_t}\bigg(uE_{k-1} -\lambda' \frac{E^2_{k-1}}{E_k}\bigg)d\mu+(n-k+1)\int_{\Sigma_t}\frac{uE_{k}-\lambda' E_{k-1}}{\Theta}~d\mu \leq 0\nonumber
\end{align}
due to the Newton inequality $E_{k} E_{k-2}\leq E_{k-1}^2$, the divergence theorem and the pinching condition \eqref{PC}.

As a consequence, the limiting hypersurface along flow \eqref{fl1}, as well as that along flow \eqref{FL-3}, is some coordinate slice, say $\Sigma_\infty:=\{r_\infty\}\times\mathbb{S}^n$, see \cite{MS}.

Along flow \eqref{FL-3}, the volume $\mathcal{A}_{-1}$ keeps constant, while the weighted volume $\mathcal{A}_{-1}$ is monotone decreasing, since
\begin{align}
\frac{d}{dt}\mathcal{B}_{-1}(\Sigma_t)=-\frac{n+1}{n}\int_{\Sigma_t}||\nabla\lambda'||^2 d\mu\leq 0\nonumber
\end{align}
due to \eqref{How} and the divergence theorem. Furthermore, along flow \eqref{FL-3}, there holds
\begin{align}
\frac{d}{dt} \mathcal{A}_l(\Sigma_t)\geq (n-l)\int_{\Sigma_t}(uE_{l+2}-\lambda'E_{l+1}) d\mu=0\nonumber
\end{align}
for $l\in\{1,\cdots, n-1\}$, where Newton's inequality and \eqref{MI} are used.

On the other hand, for each $\l\in\{1,\cdots,k-1\}$, the weighted curvature integral $\mathcal{B}_l$ is monotone decreasing along flow \eqref{fl1}, noticing that
\begin{align}
\frac{d}{dt}\mathcal{B}_l(\Sigma_t)=(1+l)\int_{\Sigma_t} uE_l\mathcal{F} d\mu+(n-l)\int_{\Sigma_t} \lambda' E_{l+1} \mathcal{F} d\mu.\nonumber
\end{align}
In the particular case where $k=n$, the weighted curvature integral $\mathcal{B}_n$ keeps constant along flow \eqref{fl1} since
\begin{align}
\frac{d}{dt} \mathcal{B}_n(\Sigma_t)=(n+1)\int_{\Sigma_t} (u E_n-\lambda' E_{n-1}) d\mu=0.\nonumber
\end{align}
On the other hand, the quermassintegral $\mathcal{A}_{n-1}$ is monotone increasing since
\begin{align}
\frac{d}{dt}\mathcal{A}_{n-1}(\Sigma_t)&=\int_{\Sigma_t} \frac{uE_{n}-\lambda' E_{n-1}}{u}~d\mu\nonumber\\
&=\frac{1}{n}\int_{\Sigma_t}\frac{E_n^{ij}(\lambda')_{;i} u_{;j} }{u^2} d\mu\nonumber
\end{align}
is non-negative due to the convexity of evolving hypersurface.

Note that in de Sitter space the coordinate slice $\{s\}\times\mathbb{S}^n$ is spacelike and convex with its principal curvatures
\begin{align}
\kappa_1(s)=\cdots=\kappa_n(s)=\frac{\lambda'(s)}{\lambda(s)}< \frac{\lambda(s)}{\lambda'(s)}.\nonumber
\end{align}
Let $\psi_{-1}(s)=\omega_n \lambda^{n+1}(s)$, and for $l\in\{1,\cdots,n\}$
\begin{align}
\psi_l(s)=\omega_n(\lambda'(s))^{l+1} \lambda^{n-l}(s) ,\nonumber
\end{align}
then $\psi_{-1}'(s)=(n+1)\omega_n \lambda^{n}(s)\lambda'(s)>0$, while for $l\in\{1,\cdots,n\}$
\begin{align}
\psi_l'(s)=\omega_n(\lambda'(s))^l \lambda^{n-l-1}(s)\big[(n+1)\lambda^2(s)-(n-l)\big]>0.\nonumber
\end{align}
Now we calculate the quermassintegrals of the coordinate slice $\{s\}\times\mathbb{S}^n$. In fact, they are equal to
\begin{align}
\varphi_{-1}(s)=(n+1)\omega_n\int_0^s \lambda^n(t)dt,~\;~ ~\;~ \varphi_0(s)=\omega_n\lambda^n(s), \nonumber
\end{align}
and
\begin{align}
\varphi_{k}(s)&=\omega_n\bigg([\lambda'(s)]^{k}\lambda^{n-k}(s)-\frac{k}{n-k+2}\varphi_{k-2}(s)\bigg), ~ k=1,\cdots,n,\nonumber
\end{align}
respectively, which are obviously invertible functions in $s$.

Consequently, there are geometric inequalities
\begin{align}
\mathcal{B}_n(\Sigma)= \psi_n(r_\infty)=\psi_n \circ\psi^{-1}_l(\psi_l(r_\infty))\leq\psi_n \circ\psi^{-1}_l(\mathcal{B}_l(\Sigma))\nonumber
\end{align}
and
\begin{align}
\mathcal{A}_{n-1}(\Sigma)\leq \varphi_{n-1}(r_\infty)=\varphi_{n-1}\circ \psi_{n}^{-1} (\psi_{n}(r_\infty))= \varphi_{n-1}\circ \psi_{n}^{-1}(\mathcal{B}_{n}(\Sigma)),\nonumber
\end{align}
where $\Sigma$ is some closed spacelike strictly convex hypersurface which satisfies \eqref{PC}, thus can be considered as the initial value of flow \eqref{fl1}.

\end{document}